\documentclass[10pt,reqno]{article}     %% Input Standard CPAM class.

\usepackage[left=1.4in, top=1in, bottom=1in, right=1.4in, asymmetric]{geometry}

\usepackage{amsmath,amsthm,amssymb,amscd,fancybox,ifthen,float,epsfig,subfigure}
\usepackage{graphicx,xcolor}
\usepackage{setspace}
\usepackage[all]{xy}
\usepackage{comment}
\usepackage{hyperref}
\usepackage{mathtools}
\usepackage{mathpazo}

\numberwithin{equation}{section}

\theoremstyle{definition}

\theoremstyle{remark}

\newcommand{\ignore}[1]{}

\floatstyle{plain}
\restylefloat{figure}

\usepackage{etoolbox}
\apptocmd{\sloppy}{\hbadness 10000\relax}{}{}

\def\Xint#1{\mathchoice
   {\XXint\displaystyle\textstyle{#1}}%
   {\XXint\textstyle\scriptstyle{#1}}%
   {\XXint\scriptstyle\scriptscriptstyle{#1}}%
   {\XXint\scriptscriptstyle\scriptscriptstyle{#1}}%
   \!\int}
\def\XXint#1#2#3{{\setbox0=\hbox{$#1{#2#3}{\int}$}
     \vcenter{\hbox{$#2#3$}}\kern-.5\wd0}}

\def\dashint{\Xint-}

\begin{document}

\title{Signal processing approach to mesh refinement in simulations of axisymmetric droplet dynamics}

\author{Kazuki Koga}

\maketitle

\begin{abstract} 
We propose a novel mesh refinement scheme based on signal processing for boundary integral simulations of inviscid droplet dynamics with axial symmetry. A key idea is to directly access the Fourier coefficients of a principal curvature as a function of the arclength through a natural change of variables. The trapezoidal rule is applied to those Fourier-type integrals and the resulting formula fits in the framework of the non-uniform fast Fourier transform. This observation enables to efficiently use an envelope analysis and smoothing filter to generate guidelines for mesh refinement in two singularity formation scenarios. Applications also include a non-iterative construction of the uniform parametrization for an important class of plane curves, which is used in a convergence study of the time-stepping procedure implemented in the previous work by Nitsche and Steen [J. Comput. Phys. 200 (2004) 299].
\end{abstract}

%%%%%%%%%%%%%%%%%%%%%%%%%%%%%%%%%%%
%%%%%%%%%%%%%%%%%%%%%%%%%%%%%%%%%%%
%%%%%%      INTRODUCTION     %%%%%%%%%%%%%%%%%
%%%%%%%%%%%%%%%%%%%%%%%%%%%%%%%%%%%
%%%%%%%%%%%%%%%%%%%%%%%%%%%%%%%%%%%

\section{Introduction}\label{sec1}
A deep understanding of droplet motion under the action of surface tension is essential in various applications such as ink-jet printing, spraying, and atomization. There are a large body of literature on theoretical, experimental, and numerical studies of the physics, and significant part of those works is devoted to investigating singularity formations observed when a fluid interface self-intersects. Such phenomena have been of great interest for their universality, richness, and elegant characterizations.\par
To reveal the structures of those singularities, a number of works have performed boundary integral simulations for potential flow models that describe the motion of an interface separating two inviscid fluids. For example, a two-dimensional droplet formation by a roll-up due to the Kelvin-Helmholtz instability is carefully studied by {Hou et al. \cite{HoLoSh1997}} (henceforth HLS97). In an axisymmetric case, Leppinen and Lister \cite{LeLi2003} investigate the self-similarity of inviscid capillary pinch-off for a range of density jumps across the surface of a droplet. Their results are more firmly confirmed by Nitsche and Steen \cite{NiSt2004} (henceforth NS04) with additional insights into the nature of the self-similarity. Furthermore, a recent work by Burton and Taborek \cite{BuTa2007B} revisits pinch-off of 2D droplets and finds that the type of the singularity is distinctly different from that of its axisymmetric counterpart.\par
Besides the dimension of the problem being reduced by one, a major advantage of boundary integral simulations is its flexibility in controlling the computational mesh. The successful works  mentioned above undoubtedly rest on mesh refinement schemes that dynamically redistribute computational points densely in the regions of interest. However, the strategies employed in the previous studies are highly problem-specific. That is, it is required to know ``a priori" a few important properties of given initial conditions, which precludes applications of those methods to general cases.\par
As an extension of NS04, this paper presents a novel mesh refinement strategy for resolving singularity formations in the vortex sheet formulation of axisymmetric droplet dynamics driven by surface tension. A vortex sheet is a sharp interface separating two inviscid, incompressible, and irrotational fluids shearing past each other, which has a boundary integral formulation that derives from the Biot-Savart law. In our scheme, guidelines for mesh refinement are generated based on an envelope analysis and smoothing filter from signal processing applied to the Fourier coefficients of the curvature in the symmetry plane as a function of the arclength. These Fourier coefficients are discretized by the trapezoidal rule after a natural change of variables, and the resulting discrete formula is approximated by the non-uniform fast Fourier transform \cite{BaMaKl2019} with $\mathcal{O}(N\log N)$ complexity. It turns out that our method is capable of detecting singularity formations with far less human intervention than the existing methods, provided that highly complex geometry does not appear in a very short time. We test the numerical scheme on the pair of initial conditions studied by Nie \cite{Nie2001} (henceforth Nie01) that reach two distinct types of singularities at the end, and discuss several aspects of our strategy related to its efficiency and accuracy. Also, applications of our ideas include a non-iterative algorithm for constructing the uniform parametrization of closed plane curves, which is utilized in a convergence study of the time-stepping procedure that has not been sufficiently justified in NS04.\par
The rest of this paper is organized as follows. Section \ref{sec2} introduces a boundary integral formulation of axisymmetric vortex sheets with surface tension and two relevant initial conditions. Section \ref{sec3} describes the theory behind our mesh refinement strategy based on signal processing. Section \ref{sec4} explains our numerical method and some technicalities in its implementation. Section \ref{sec5} presents numerical results and validates the implemented scheme. Concluding remarks are given in Section \ref{sec6}.
%utility, prototypical, encoded, by virtue of, albeit with, nonetheless, ambient, extrapolated to find $t_p$,

%%%%%%%%%%%%%%%%%%%%%%%%%%%%%%%%%%%
%%%%%%%%%%%%%%%%%%%%%%%%%%%%%%%%%%%
%%%%%%      FORMULATION       %%%%%%%%%%%%%%%%%
%%%%%%%%%%%%%%%%%%%%%%%%%%%%%%%%%%%
%%%%%%%%%%%%%%%%%%%%%%%%%%%%%%%%%%%

\section{Formulations}\label{sec2}
In this section, we give a brief introduction to boundary integral equations for axisymmetric vortex sheets with surface tension, which is followed by a short description of finite-time singularities in the formulation reported by other authors. For interested readers, a more detailed derivation can be found in NS04.\par
\subsection{Lagrangian theory}
We consider irrotational motion of two inviscid, incompressible, and immiscible fluids in $\mathbb{R}^3$
separated by a closed surface $S$. In the following, the inner and outer fluid are denoted by $\Omega_1$ and $\Omega_2$, respectively, and $\Omega_1$ is assumed to be simply-connected. For simplicity, we also assume that the density of each fluid is uniformly equal to unity and that no external force acts on the system. Then, the fluid velocity $\mathbf{u}_i$ and pressure $p_i$ in $\Omega_i$ ($i=1,2$) satisfy the Euler equation with the boundary conditions 
%the divergence-free condition, and the no-vorticity condition :
%
%\begin{equation}
%\label{eq:def_Euler}
%\frac{\partial \mathbf{u}_i}{\partial t} + (\mathbf{u}_i \cdot \nabla)\mathbf{u}_i+\nabla p_i=0, \quad \nabla \cdot \mathbf{u}_i=0, \quad \nabla \times \mathbf{u}_i=0,\quad \mbox{in } \Omega_i
%\end{equation} 
%
%\noindent where $t$ is the time variable. Besides, we add the following boundary conditions that determine $\mathbf{u}_i$ %and characterize the motion of $S$ :
%
\begin{alignat}{3}
\label{eq:def_bd_farfield}&\mathbf{u}_i(\mathbf{x},t) \rightarrow 0\;\; \mbox{ as }\; \|\mathbf{x}\| \rightarrow \infty \quad \mbox{(far-field condition),}\\[5pt]
\label{eq:def_bd_kinem}&[\mathbf{u}]_S \cdot\mathbf{n}=0 \quad \mbox{on } S\quad \mbox{(kinematic boundary condition)},\\[5pt]
\label{eq:def_bd_LaYo}&[p]_S =\sigma \kappa \quad\;\;\; \mbox{on } S\quad \mbox{(Laplace-Young condition)}.
\end{alignat}
where $\mathbf{x}$ and $t$ are the spatial and time variable, respectively. Also, $\mathbf{n}$ is the unit normal vector of $S$ pointing to $\Omega_1$, $\kappa$ is twice the mean curvature that is positive if $S$ is a sphere, and $\sigma$ is a positive surface tension coefficient. The brackets $[\cdot]_S$ measure the jump of a function across $S$ in the following manner: 
\begin{equation}
\label{eq:def_bracket}
[f]_S(\mathbf{x},t) =\lim_{\substack{\mathbf{y}\rightarrow \mathbf{x} \\ \mathbf{y}\in \Omega_1}}f_1(\mathbf{y},t) -\lim_{\substack{\mathbf{y}\rightarrow \mathbf{x} \\ \mathbf{y}\in \Omega_2}}f_2(\mathbf{y},t), \quad \mathbf{x}\in S.
\end{equation}
The far-field condition (\ref{eq:def_bd_farfield}) states that the fluid at infinity is at rest. The kinematic boundary condition (\ref{eq:def_bd_kinem}) requires that the normal component of the fluid velocity must be continuous across $S$, while the tangential component can be discontinuous there. The Laplace-Young condition (\ref{eq:def_bd_LaYo}) represents the effect of surface tension on $S$ in terms of the mean curvature and the pressure jump across the interface.\par
\begin{figure}[t]
\begin{center}
\includegraphics[width=\linewidth,trim=0 0 0 0]{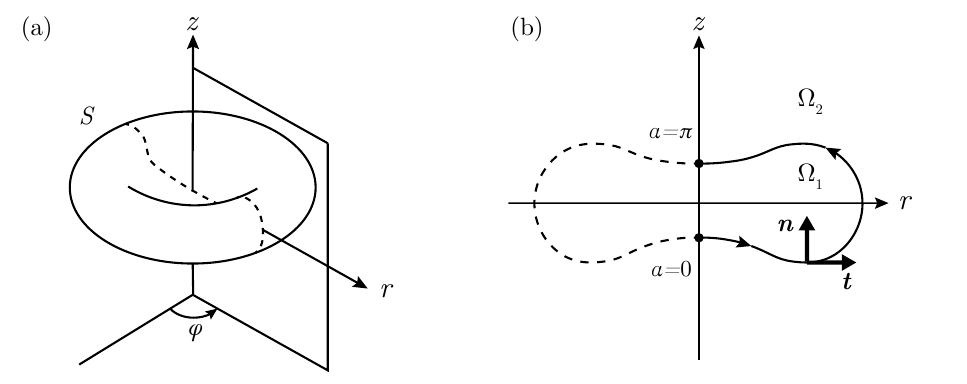}
\end{center}
\caption{\label{fig:axisym} Sketch illustrating: (a) axial symmetry, (b) parametrization in $r$--$z$ plane.}
\end{figure}
Throughout the rest of this paper, we assume that the flow is axially symmetric. That is, in some cylindrical coordinates $(r,\varphi,z)$, all functions describing the flow are independent of the azimuthal coordinate $\varphi$ (see Fig.\ref{fig:axisym}(a)). In this setting, for tracking the motion of a free surface $S$, it suffices to solve time evolution of the intersection of $S$ and a symmetry plane $\varphi =\mbox{const.}$ Thus, we parametrize the curve in the $r$--$z$ plane by
%
%%---Eq. Parametrization in (r,z)---%%
\begin{equation}
\label{eq:def_prmt}
\mathbf{X}(\alpha,t) = (r(\alpha,t), z(\alpha,t)),\quad   \alpha \in [0,\pi],
\end{equation}
where $r$ and $z$ are radial and axial coordinates, respectively, and find its evolution equation
%
%%---Evolution Eq. for X---%%
\begin{equation}
\label{eq:decompX}
\frac{\partial \mathbf{X}}{\partial t} = U\mathbf{n}+V\mathbf{t}.
\end{equation}
Again, $\mathbf{n}$ is the unit normal vector pointing to $\Omega_1$, and $\mathbf{t}$ is the unit tangent vector in the counterclockwise direction. In particular, we set $r(0,t)=r(\pi,t)=0$ for all $t$, and direct $\mathbf{X}$ so that it satisfies $\mathbf{t}=\mathbf{X}_\alpha/s_\alpha$ (see Fig.\ref{fig:axisym}(b)). Here, the subscript $\alpha$ denotes differentiation with respect to the variable $\alpha$, and $s_\alpha =\sqrt{r_\alpha^2 + z_\alpha^2}$ is called the {\it relative spacing} that is the derivative of the arclength $s$ measured from the starting point $\mathbf{X}(0,t)$.  In the present case, the kinematic boundary condition (\ref{eq:def_bd_kinem}) implies that $U=\mathbf{W} \cdot \mathbf{n}$, where $\mathbf{W}$ is the average of the fluid velocities on either side of $S$, whereas it puts no restriction on $V$.  In fact, it is well-known that the shape of an evolving plane curve is determined solely by the normal velocity $U$, and the tangential velocity $V$ only affects the relative spacing $s_\alpha$. \par
To obtain an integral form of $\mathbf{W}$, we introduce a dynamical variable $\gamma$ called the {\it unnormalized vortex sheet strength} that measures the jump of the tangential component of the fluid velocity in a frame-dependent manner:
\begin{equation}
\label{eq:def_gamma}
\gamma(\alpha,t) =\mathbf{X}_\alpha \cdot [\mathbf{u}]_S(\mathbf{X}(\alpha,t), t),\quad  \alpha \in [0,\pi].
\end{equation}
Then, assuming that the azimuthal component of the fluid velocity is uniformly zero ({i.e. without swirl}), the Biot-Savart law and the Plemelj's formula imply that $\mathbf{W}=(w_r,w_z)$ is given by the principal-values (denoted by dashed integrals)
%%
%%---Biot-Savart in (r,z)---%%
\begin{alignat}{2}
\label{eq:def_axbiot_r}
&w_r(\alpha,t)= \frac{1}{2\pi}\dashint_0^{\pi}\biggl(\frac{\gamma'}{\rho_2}\biggr) \biggl(\frac{z'-z}{r}\biggr)\biggl[K(\lambda)-\frac{(z'-z)^2+r^2+r'^2}{\rho^2_1}E(\lambda)\biggr] d\alpha',\\
\label{eq:def_axbiot_z}
&w_z(\alpha,t)= \frac{1}{2\pi}\dashint_0^{\pi}\biggl(\frac{\gamma'}{\rho_2}\biggr)\biggl[K(\lambda)-\frac{(z'-z)^2+r^2-r'^2}{\rho^2_1}E(\lambda)\biggr]d\alpha',
\end{alignat}
where, for example, $r=r(\alpha,t)$ and $r'=r(\alpha',t)$ (e.g., see \cite{CaLi1992} for details). Here, $K$ and $E$ are the complete elliptic integrals of the first and second kind, respectively. Other functions in (\ref{eq:def_axbiot_r}) and (\ref{eq:def_axbiot_z}) are defined as
%
%%---Components of Biot-Savart in (r,z)---%%
\begin{equation}
\label{eq:axbiot_comp}
 \rho_1^2= (z'-z)^2+(r'-r)^2,\quad \rho_2^2= (z'-z)^2+(r'+r)^2, \quad \lambda^2 ={4rr'}/{\rho^2_2}.   
\end{equation}
On the other hand, the evolution equation for $\gamma$ is a consequence of the unsteady Bernoulli's theorem and the Laplace-Young condition (\ref{eq:def_bd_LaYo}):
%
%%---Evolution Eq. for Strength---%%
\begin{equation}
\label{eq:dt_gamma}
\frac{\partial \gamma}{\partial t} = -\sigma \kappa_\alpha + ((V-\mathbf{W}\cdot \mathbf{t})\gamma/s_\alpha)_\alpha.
\end{equation}
In axisymmetric geometry, it is convenient to express $\kappa$ as the sum of
%
%%---Def. of Curvatures---%%
\begin{equation}
\label{eq:curvature}
\kappa_z=\frac{z_{\alpha\alpha}r_{\alpha}-z_{\alpha}r_{\alpha\alpha}}{s_\alpha^3} ,\quad \kappa_r = \frac{z_\alpha}{s_\alpha r},
\end{equation}
which are the principal curvatures in the symmetry plane and in a plane normal to it, respectively. Without surface tension (i.e. $\sigma=0$), the first term in the right-hand side of (\ref{eq:dt_gamma}) is absent, and choosing $V=\mathbf{W}\cdot \mathbf{t}$ renders $\gamma$ time-independent. This special choice for $V$ is often referred to as the {\it Lagrangian frame}, and it has been studied mainly in the context of the famous Moore's singularity \cite{Moore1979}. With surface tension, however, it is no longer straightforward to keep $\gamma$ conserved in time. Moreover, it is known that the Lagrangian frame typically leads to rapid decreases in $s_\alpha$ at points where geometry is not necessarily complex, as demonstrated by HLS97 and NS04.  These facts motivate us to choose $V$ to control the relative spacing $s_\alpha$ for other practical purposes. \par
%
%%%%%%      ALTERNATIVE VARIABLES         %%%%%%%%%%%%%%%%%
%
Following {Hou et al. \cite{HoLoSh1994}} (henceforth HLS94), we switch from the Cartesian coordinates $(r,z)$ to the so-called angle--arclength variables $(\theta,s_\alpha)$ defined via the relation
%
%%---Evolution Eq. for Theta---%%
\begin{equation}
\label{eq:def_theta}
\mathbf{X}_\alpha =(s_\alpha\cos\theta, s_\alpha \sin \theta).
\end{equation}
Having these new dynamical variables, we are able to reconstruct $(r,z)$ up to a translation by integrating Eq.(\ref{eq:def_theta}). Such a loss of information is not a problem here, for the current interest lies in developing mesh refinement schemes to resolve singularity formations due to self-intersections of the interface $S$. The evolution equations for $(\theta,s_\alpha)$ follow from the orthogonal decomposition in (\ref{eq:decompX}), the Frenet-Serret formula for plane curves, and the identity $\theta_\alpha =s_\alpha \kappa_z$:
%
%%---Evolution Eq. for Theta---%%
\begin{alignat}{2}
\label{eq:ev_theta}
\frac{\partial\theta}{\partial t}=\frac{U_\alpha+V\theta_\alpha}{s_\alpha}, \\
\label{eq:evol_sa}
\frac{\partial s_\alpha}{\partial t}=V_\alpha-\theta_\alpha U.
\end{alignat}
Eq.(\ref{eq:evol_sa}) is central to the present work. As mentioned earlier, because of the tangential velocity being arbitrary, the relative spacing $s_\alpha$ can be controlled by a user-specified $V$ while keeping the image of $\mathbf{X}$ unchanged. In fact, this can be done by deriving a desirable $s_\alpha$ first and calculating $V$ accordingly based on Eq.(\ref{eq:evol_sa}). Our choice for $V$ is described in detail in Section \ref{sec3}.
%
%An advantage of the variables $(\theta,s_\alpha)$ is that the relation $\theta_\alpha =s_\alpha \kappa_z$ enables to approximate the curvature $\kappa_z$, which plays a central role in the later sections, as the first derivative of $\theta$, rather than higher derivatives of $(r,z)$. Furthermore, it allows to group out the most singular term hidden in $U$ as a linear term and treat it implicitly \cite{HoLoSh1994}. Unfortunately, we do not benefit from the second advantage because of the time-stepping procedure described later, and use the new variables simply for the sake of high accuracy.
%
%%%%%%      SINGULARITY FORMATION         %%%%%%%%%%%%%%%%%
%
\subsection{Singularity formations}
It is well-known that, without any regularization, the dynamics of vortex sheets is subject to the Kelvin-Helmholtz instability and generally leads to the Moore's singularity \cite{Moore1979}. On the other hand, HLS94 points out that surface tension acts as a dispersive regularization of the Kelvin-Helmholtz instability, and also reports that numerical solutions with $\sigma\neq0$ continue beyond the Moore's singularity for the initial condition studied by Krasny \cite{Krasny1986} in the zero surface tension setting. Instead of the Moore's singularity, however, several works have found that other kinds of singularities appear when a fluid interface self-intersects. Such phenomena are often called {\it topological singularities}. \par 
%
%
%%---Fig. Singularity Formation and Curvature Profile---%%
\begin{figure}[t]
\begin{center}
\includegraphics[width=\linewidth,trim=0 0 0 0]{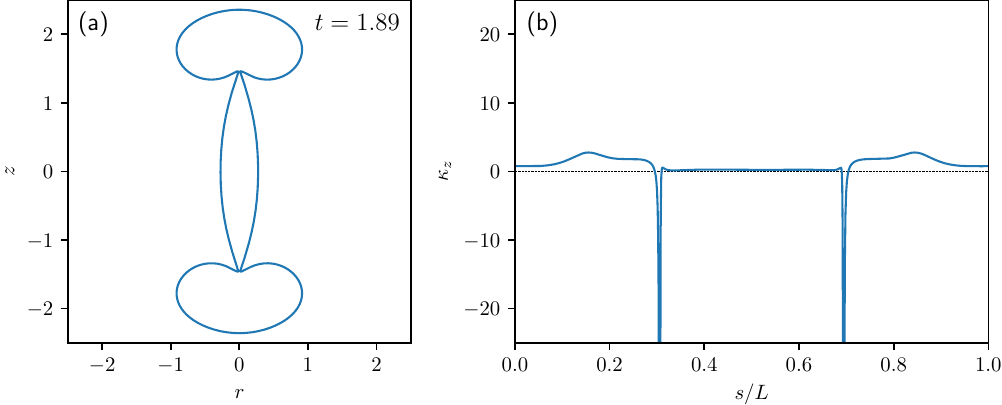}
\end{center}
\caption{\label{fig:slice_kappa_pinchoff} Solution close to singularity for pinch-off problem (\ref{eq:init_nie2}): (a) interface position in $r$--$z$ plane, (b) profile of curvature $\kappa_z$.}
\end{figure}
The most famous topological singularity is presumably {\it capillary pinch-off}, which is characterized by a fast drop in the radius of a neck due to surface tension acting on the azimuthal curvature $\kappa_r$. To study this process numerically with the vortex sheet formulation, Nie01 proposes the initial condition 
%
%%---Initial Condition for Symmetric Pinchoff---%%
\begin{equation}
\label{eq:init_nie2}
\sigma = 0.2, \quad \theta(\alpha,0) = \alpha, \quad  s_\alpha(\alpha,0)=1, \quad \gamma(\alpha,0) = -2\sin(2\alpha).
\end{equation}
Qualitatively, this problem resembles an experiment by Robinson and Steen \cite{RoSt2001}. As an intuitive illustration, we show in Fig.\ref{fig:slice_kappa_pinchoff} a snapshot of a numerical solution starting from (\ref{eq:init_nie2}) and the corresponding curvature $\kappa_z$ as a function of $s$ divided by the total arclength $L=s(\pi,t)$. As seen in Fig.\ref{fig:slice_kappa_pinchoff}(a), the solution is symmetric with respect to the $r$--axis and eventually forms two end-droplets at the top and bottom along with an elongated satellite droplet in between. Also, one can easily see that overturning occurs in each pinch-off region, which precludes descriptions of the process by a single-valued function of $z$. NS04 studies the same problem numerically to confirm the self-similarity of inviscid capillary pinch-off
\begin{equation}
\label{eq:scaling_nie2}
r_\text{min} \sim (t_p-t)^{\frac{2}{3}},\quad(z_\text{min}-z_p) \sim (t_p-t)^{\frac{2}{3}},
\end{equation}
 for a range of density jumps across the interface. Here, $(r_\text{min},z_\text{min})$ are the coordinates of a point that attains the minimum neck radius, $z_p$ is the axial coordinate where the pinch-off occurs, and  $t_p$ is the time when the radius $r_\text{min}$ reaches exactly zero. From a numerical simulation perspective, the initial condition (\ref{eq:init_nie2}) is particularly suitable for an accurate numerical study on the scaling law (\ref{eq:scaling_nie2}), for it has an additional line symmetry and the shape of the interface is considerably simple away from the necks (see Fig.\ref{fig:slice_kappa_pinchoff}(b)). In this paper, this problem is revisited to assess efficiency and accuracy of our numerical schemes.
\par
\begin{figure}[t]
\begin{center}
\includegraphics[width=\linewidth,trim=0 0 0 0]{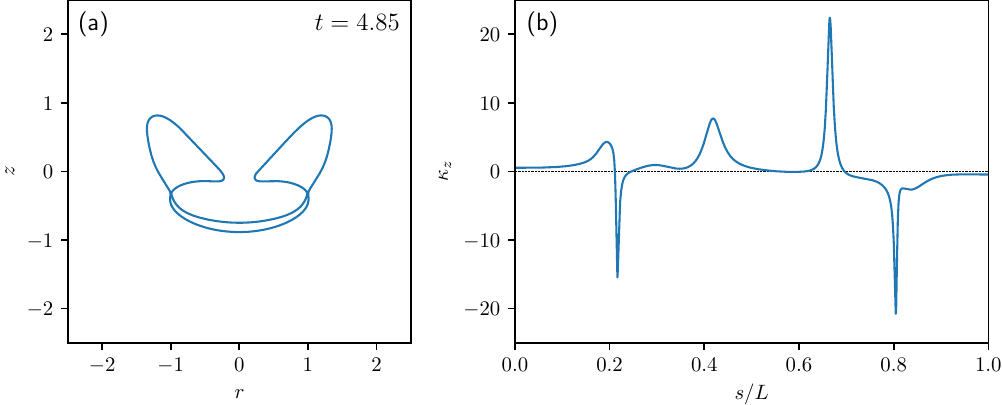}
\end{center}
\caption{\label{fig:slice_kappa_inertia} Solution close to singularity for axisymmetric bag breakup (\ref{eq:init_nie1}): (a) interface position in $r$--$z$ plane, (b) profile of curvature $\kappa_z$.}
\end{figure}

Another scenario of topological singularity formations is self-intersections due to inertia resisted by surface tension. In the present formulation, this singularity is found for the initial condition 
%
%%---Initial Condition for Bag Breakup---%%
\begin{equation}
\label{eq:init_nie1}
\sigma = 0.04, \quad \theta(\alpha,0) = \alpha, \quad  s_\alpha(\alpha,0)=1, \quad \gamma(\alpha,0) = -\sin(\alpha),
\end{equation}
which is also suggested by Nie01. This situation can be understood as a uniform flow past a sphere \cite{Batchelor1967} where the solid boundary is instantly dissolved at $t=0$ and evolves under the action of surface tension. Again, Fig.\ref{fig:slice_kappa_inertia} plots a numerical solution starting from (\ref{eq:init_nie1}) and its curvature $\kappa_z$. 
%
%
%%---Fig. Singularity Formation and Curvature Profile---%%
As seen in Fig.\ref{fig:slice_kappa_inertia}(a), the interface collides itself away from the axis of symmetry, and the droplet seems to break into a torus-like component and a simply-connected one. For brevity, we refer to this type of singularities as {\it bag breakup}. The nature of bag breakup is expected to be different from that of capillary pinch-off, because in this case surface tension should always tend to restore the droplet to a spherical shape. Regarding this point, motivated by their own experiments \cite{BuTa2007A}, Burton and Taborek \cite{BuTa2007B} numerically investigate bag breakup in two-dimensional flows (called {\it 2D pinch-off} by the same authors), and they find that its self-similarity is of the second kind. That is, the local behavior of a fluid interface close to a self-intersection is characterized by the memory of an initial condition, and the spatial variables follow scaling laws with different exponents. Furthermore, by using the full model and its ``slender" approximation, they find numerical evidences of those self-similarity exponents distinct from that of (\ref{eq:scaling_nie2}).  We speculate that a similar analysis may be applicable to the axisymmetric case. In this paper, however, we do not aim to identify the similarity exponents in the axisymmetric bag breakup (\ref{eq:init_nie1}), because it forms multiple high-curvature regions away from the self-intersection and is therefore inappropriate for that purpose (see Fig.\ref{fig:slice_kappa_inertia}(b)). Rather, our interest lies in whether it is possible to develop a mesh refinement scheme that is able to capture singularity formations and such high-curvature regions simultaneously.
%%%%%%%%%%%%%%%%%%%%%%% %%%%%%%%%%%%
%%%%%%      REFINEMENT          %%%%%%%%%%%%%%%%%
%%%%%%%%%%%%%%%%%%%%%%% %%%%%%%%%%%%

\section{Mesh refinement}\label{sec3}
The choice of the tangential velocity $V$, which we have left undefined so far, is crucial in numerical simulations of evolving plane curves. This is particularly true if one needs to resolve small-scale structures of a fluid interface approaching an imminent singularity formation. In the context of potential flow models, several works use the Lagrangian frame $V=\mathbf{W}\cdot \mathbf{t}$ in a combination with some regridding techniques. For example, in simulations of capillary pinch-off with axial symmetry, Leppinen and Lister \cite{LeLi2003} advect computational points by the Lagrangian frame and redistribute them at each time step in order to keep the relative spacing proportional to the distance from the point of the minimum neck radius. A similar strategy is used by Burton and Taborek \cite{BuTa2007B} to study bag breakup in two-dimensional flows. Although this idea provides a simple way of accumulating points to high-curvature regions, it does not naturally extend to general cases because their methods heavily rely on the fact that the initial conditions studied in those works lead to self-intersections right on a targeted coordinate axis.  \par 
In the following, we employ a style in mesh refinement that harnesses the arbitrariness of the tangential velocity $V$. Namely, the equation for $V$ that realizes a given $s_\alpha$ is derived for a simple open curve identified with an axisymmetric surface whose inner region is simply connected, which implies that
\begin{equation}
\label{eq:def_Vend}
V(0,t)=V(\pi,t)=0,
\end{equation}
for all $t$. At a glance, this condition might seem to limit the applicability of our method. However, with minor modifications to (\ref{eq:def_Vend}), we believe that our ideas easily generalize to the cases of periodic curves and unbounded curves with periodic boundary conditions. \par
 Following HLS94, we define the relative spacing $s_\alpha$ as a product of the total arclength $L$ and a strictly positive function $R$ whose integral over $[0,\pi]$ is equal to 1:
 %%---Def. of s_\alpha with Ratio Function.---%%
\begin{equation}
\label{eq:ratio}
s_\alpha(\alpha,t)=R(\alpha,t)L(t), \quad \int_0^{\pi} Rd\alpha=1.
\end{equation} 
By differentiating both sides of the first equation in (\ref{eq:ratio}) with respect to $t$ and substituting the result into (\ref{eq:evol_sa}), we get 
%
%%---Differential Equation for V---%%
\begin{equation}
\label{eq:diffeq_V}
V_\alpha= \biggl(\frac{\partial R}{\partial t}L+ R\frac{dL}{d t}\biggr)+ \theta_\alpha U.
\end{equation} 
Also, by combining (\ref{eq:evol_sa}) and (\ref{eq:def_Vend}), it is easy to see that
%
%%---Differential Equation for V---%%
\begin{equation}
\label{eq:evol_L}
\frac{dL}{dt}= -\int_0^{\pi}\theta_\alpha U d\alpha.
\end{equation} 
Now, by integrating (\ref{eq:diffeq_V}) with (\ref{eq:def_Vend}), we obtain the following integral equation for $V$ :
%
%%---Integral Equation for V---%%
\begin{equation}
\label{eq:inteq_V}
V(\alpha,t)= \int_0^{\alpha}\biggl(\frac{\partial R'}{\partial t}L+ R'\frac{dL}{d t}\biggr)d\alpha'+ \int_0^{\alpha}\theta_\alpha 'U'd\alpha'.
\end{equation} 
The uniform parametrization, by which we mean that $s_\alpha$ is independent of the variable $\alpha$, corresponds to $R\equiv \pi^{-1}$ \cite{HoLoSh1994}. In this case, Eq.(\ref{eq:inteq_V}) is explicitly solved for $V$ as
%
%%---T for Uniform Paramet.---%%
\begin{equation}
\label{eq:exp_uniV}
V(\alpha,t)=-\frac{\alpha}{\pi}\int_0^{\pi}\theta_\alpha Ud\alpha+\int_{0}^{\alpha}\theta_\alpha' U'd\alpha'.
\end{equation}
The uniform parametrization is a typical choice for various moving boundary problems because, unlike the Lagrangian frame, it prevents computational points from clustering at undesired locations and maintains stable spatial discretization. However, it is almost obvious that this is not the best option for accurate numerical studies of singularity formations in general. From this point of view, for a two-dimensional droplet formation by a Kelvin-Helmholtz roll-up, HLS97 alternatively proposes a variable but time-independent $R$ with a few minima, which is empirically constructed from numerical results for the same problem with the uniform parametrization. A major advantage of this method is that it is still possible to solve Eq.(\ref{eq:inteq_V}) explicitly for $V$ in the same form as (\ref{eq:exp_uniV}), but such a function $R$ is literally problem-specific by definition. \par
As an extension of HLS97, NS04 suggests a strategy for constructing a variable and time-dependent $R$ to solve the pinch-off problem (\ref{eq:init_nie2}). Their choice for $R$ is initially uniform and develops exactly two symmetric minima at $\alpha=\alpha_c$ and $\alpha=\pi-\alpha_c$. More specifically, they first prescribe the function
\begin{equation}
\label{eq:def_NiSt_guide}
g(x)=1.125+\delta(\cos(4x)-\cos(2x)),\quad \delta=\delta_\text{max}\frac{t}{t_p},
\end{equation}
for some constant $\delta_\text{max}$, and compose it with a monotonically increasing function $q(s)$ that also evolves in time and allows the two minima of $g(q)$ to dynamically follow the points of the minimal neck radius.  Using this composite function $g(q)$, they define the function $R$ as
\begin{equation}
\label{eq:def_ratio_gp}
R(\alpha,t)=\frac{g(q(s(\alpha)))}{\int_0^\pi g(q(s(\alpha'))) d\alpha'},
\end{equation}
and apply it to the initial condition (\ref{eq:init_nie2}) for a range of density jumps across the interface. Although their computations are successful in many aspects, it has to be said that this definition is still specialized to a class of initial conditions with a line symmetry. Moreover, it requires to know in advance the pinch-off time $t_p$, the number of singular points, and even a characterization of the point $(r_\text{min}, z_\text{min})$ that never causes $q(s)$ to be discontinuous in time.
  \par
The main difficulty in the design of mesh refinement for singular interfacial problems is that the curvature $\kappa_z$, a pointwise measure of the complexity of a plane curve, is by itself unlikely to define a promising candidate for $R$. Intuitively, there are two reasons for this consideration. Firstly, as shown in Fig.\ref{fig:slice_kappa_pinchoff}(b), a peak in $\kappa_z$ that corresponds to a self-intersection can be significantly sharp and therefore, if $R$ is written in terms of $\kappa_z$, the relative spacing $s_\alpha$ may reflect this sharpness as a dynamical variable. Secondly, it is found in Fig.\ref{fig:slice_kappa_inertia}(b) that there is an inflection point (i.e. $\kappa_z=0$) in a small interval where the gradient of $\kappa_z$ is very steep, which causes $s_\alpha$ to increase locally and magnifies $\kappa_{z,\alpha}$ many times. A similar situation can be seen in Fig.\ref{fig:slice_kappa_pinchoff}(b), and it is numerically verified in NS04 that there exists at least one inflection point in the self-similar regime (\ref{eq:scaling_nie2}). In fact, to our knowledge, there is no previous study that succeeds in resolving self-intersections of a fluid interface with mesh refinement based on pointwise evaluations of $\kappa_z$. Nevertheless, it is no doubt that the graphs of the curvature $\kappa_z$ visually show where singularities are about or at least likely to appear. In some sense, the previous works draw rough sketches of such visual cues as their function $R$ with the help of numerical experiments or theoretical considerations.   \par 
We attempt to further extend the idea of NS04 by developing a new tool that automatically generates guideline functions, denoted by $GL$, in place of the composite function $g(q)$. In particular, we automate the ``sketching" process by applying an envelope analysis and smoothing filter from signal processing to the curvature $\kappa_z$ as a periodic function of the arclength $s$. For this purpose, the first step is to extend the parametric curve $\mathbf{X}$ to the interval $[0,2\pi]$ with the symmetry
\begin{equation}
r(\pi+\alpha,t)=-r(\pi-\alpha,t),\quad z(\pi+\alpha,t)=z(\pi-\alpha,t),
\end{equation}
as illustrated in Fig.\ref{fig:axisym}(b), and consider the Fourier coefficients
 %%---Fourier Coefficients in Arclength---%%
\begin{equation}
\label{eq:def_ft_arc}
\hat{f}(k)=\frac{1}{L_p}\int_0^{L_p}f(s)e^{-2\pi i k\frac{s}{L_p}}ds,
\end{equation}
where $L_p$ denotes the total arclength of the periodic curve $\mathbf{X}$. This quantity can be accessed via the natural change of variables $ds = s_\alpha d\alpha$ :
 %%---Fourier Coefficients in Arclength---%%
\begin{equation}
\label{eq:changevar}
\frac{1}{L_p}\int_0^{L_p}f(s)e^{-2\pi i k\frac{s}{L_p}}ds = \frac{1}{L_p}\int_0^{2\pi}(f(s(\alpha))s_\alpha)e^{-2\pi ik \frac{s(\alpha)}{L_p}}d\alpha,
\end{equation} 
and is expected to be invariant under a change of parametrization that preserves the starting point $\mathbf{X}(0,t)$ and the direction of increasing $\alpha$. This property is essential as an element of mesh refinement, because we cannot construct an appropriate $R$ from any frame-dependent information.  \par
\begin{figure}[t]
\begin{center}
\includegraphics[width=\linewidth,trim=0 0 0 0]{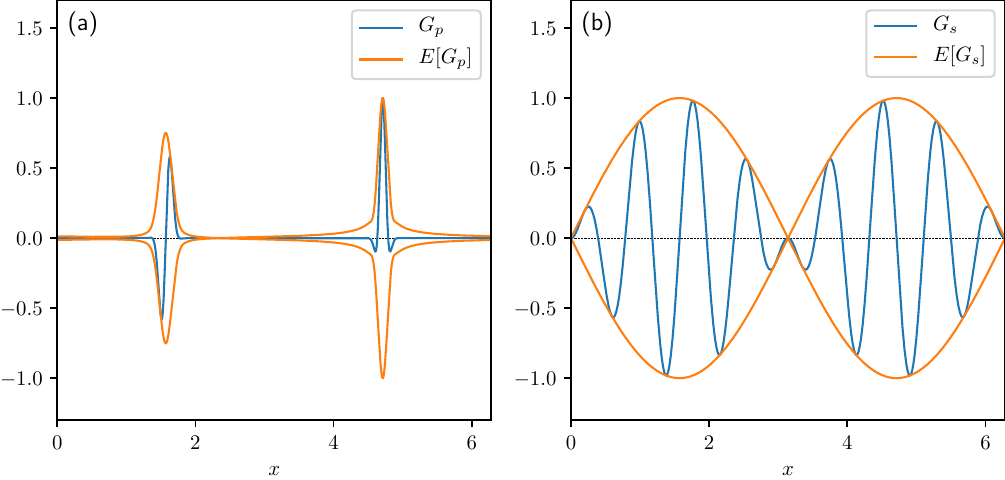}
\end{center}
\caption{\label{fig:anal_env} Signals and corresponding analytic envelopes: (a) Gaussian pulses $G_p$, (b) product of two sine functions $G_s$.}
\end{figure}

Next, we introduce the so-called {\it analytic signal} \cite{Gabor1946}
 %%---Def. of Analyticl Signal---%%
\begin{equation}
\label{eq:def_asig}
A[f](s)=f(s)+i\mathcal{H}[f](s),
\end{equation}
where $\mathcal{H}$ is the Hilbert transform, which is the Fourier multiplier of the form
 %%---Def. of Hilbert Transform---%%
\begin{equation}
\label{eq:hilbert}
\widehat{\mathcal{H}[f]}(k) = (-i \mbox{sgn}(k))\cdot \hat{f}(k). 
\end{equation}
Hence, the Fourier coefficients of $A[f]$ vanish for all negative wavenumbers, and the term ``analytic" derives from the fact that such a periodic function, if it is considered on the unit circle, can be extended to a holomorphic function on the unit disk. In classical harmonic analysis, the operator $A$ plays important roles in the study of the convergence of the Fourier series for $p$-integrable functions with $1<p<\infty$ \cite{MuSc2013}. On the other hand, in signal processing, the analytic signal is often called a {\it pre-envelope}, because, by rewriting $A[f]$ in the polar form
 %%---Def. of Analytical Envelope---%%
\begin{equation}
\label{eq:asig_phase}
A[f](s)=E[f](s) e^{i\psi(s)},
\end{equation}
one can extract the {\it analytic envelope} \cite{Rice1982}
 %%---Def. of Analytical Envelope---%%
\begin{equation}
\label{eq:aenv}
E[f](s)=\sqrt{f(s)^2+(\mathcal{H}[f](s))^2}.
\end{equation}
Here, the angle $\psi$ is called the instantaneous phase, which is also defined in terms of $f$ and its Hilbert transform. Roughly speaking, the analytic envelope is expected to have a low frequency content in comparison to that of the complex exponential factor in (\ref{eq:asig_phase}). We explain this property with the simplest example presented by Cohen \cite{Cohen1995}. Suppose that we have a smooth complex-valued function $A(x)=E(x)e^{ik_0 x}$ with a $2\pi$-periodic real-valued $E$ and an integer $k_0>0$. Then, the Fourier coefficients of $A$ are those of ${E}$ shifted to the right by $k_0$, and, for $A$ to be ``analytic", the coefficients of $E$ must be zero for wavenumbers $k< -k_0$. By the conjugate symmetry for real-valued functions, this means that all the non-zero Fourier coefficients of $E$ are contained within the range $-k_0\leq k\leq k_0$, and the strict inequalities hold if the function $A$ is zero-mean. Unfortunately, it is very hard to show the same property for periodic functions in general. Instead, to illustrate some tendencies, the analytic envelope is demonstrated in Fig.\ref{fig:anal_env} for the following functions:
\begin{alignat}{2}
&G_p(x)=e^{-64(x-C_1)^2}\sin(\omega(x-C_1))+e^{-64(x-C_2)^2}\cos(\omega(x-C_2)),\\
&G_s(s)=\sin(x)\sin(8x),
\end{alignat}
where $C_1=\frac{\pi}{2}$, $C_2=\frac{3}{2}\pi$, and $\omega=16$. 
The signal $G_p$ is a superposition of two Gaussian pulses \cite{BaFrFrKeSh1984} that have compact numerical supports because of the fast decaying prefactors. We show the raw $G_p$ and its analytic envelope in Fig.\ref{fig:anal_env}(a). Here, note that the left pulse has one zero point between two peaks with different signs, while the right pulse has one sharp peak surrounded by multiple zero points. In both cases, the operator $E$ ignores all these zero points and transforms the original pulses to simple and strictly positive peaks. Therefore, it is expected that the analytic envelope can remove the difficulty caused by inflection points near self-intersections. On the other hand, the signal $G_s$ is a product of sine functions with a high and low frequency, as shown in Fig.\ref{fig:anal_env}(b). In this case, the operator $E$ fails to detect any local complexity, presumably due to the low-frequency nature of the analytic envelope, and simply returns a function similar to $|\sin x|$. These results imply that the analytic envelope is effective only if high-curvature regions are localized and sufficiently isolated from each other. Moreover, the second example shows that the analytic envelope can be non-differentiable at zero points of the original signal. To avoid this problem, we use a regularized envelope
 %%---Def. of Regularized Analytical Envelope---%%
\begin{equation}
\label{eq:aenv_reg}
E_r[f](s)=\sqrt{1+(E[f](s))^2}.
\end{equation}
Adding a constant inside the square root renders the final form of $GL$ slightly less sharp, but this regularization is clearly inevitable for the sake of smoothness. Also, it should be noted that the analytic envelope of a nearly-singular function is inherently very sharp. In this paper, to obtain $GL$ smoother than the original $\kappa_z$, we apply the Gaussian filter to the analytic envelope $E_r[\kappa_z]$. That is, the function $E_r[\kappa_z]$ is convolved with the heat kernel 
\begin{equation}
\label{eq:def_heat}
H_a(s)=(a^2/\pi)^{\frac{1}{2}}e^{-a^2 s^2}, 
\end{equation}
where the parameter $a$ controls the level of smoothing and essentially determines the locality of mesh refinement. To summarize, we suggest the following guideline function $GL$ as a ``rough" measure of the complexity of an evolving plane curve:
 %%---Def. of Guideline Function---%%
\begin{equation}
\label{eq:guideline}
GL(s)=(E_r[\kappa_z]*H_a)(s)
\end{equation}
%
%%---Fig. Snapshots of Capillary Pinchoff---%%
\begin{figure}[t]
\begin{center}
\includegraphics[width=\linewidth,trim=0 0 0 0]{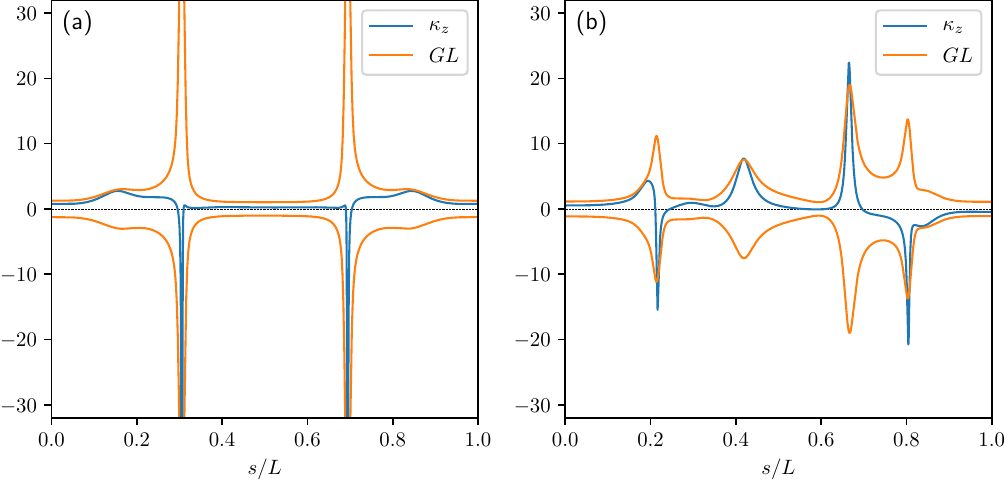}
\end{center}
\caption{\label{fig:kappa_GL} Curvature $\kappa_z$ and generated guideline function $GL$ at times close to singularity: (a) pinch-off problem at $t=1.89$, (b) axisymmetric bag breakup at $t=4.85$.}
\end{figure}
In Fig.\ref{fig:kappa_GL}, the functions $GL$ with $a=20$ are shown for the curvature $\kappa_z$ in Fig.\ref{fig:slice_kappa_pinchoff}(b) and Fig.\ref{fig:slice_kappa_inertia}(b). In both cases, our technique is successful in generating simple and strictly positive functions that outline relevant features of the original signals $\kappa_z$. \par
Finally, we define the relative spacing $s_\alpha$ in terms of the guideline function $GL$. In our definition, the function $R$ is split into a few components using a parameter $\delta_R$, which determines the maximal amount of refinement, and another Gaussian as a function of the time $t$:
 %%---Def. of Ratio Function---%%
\begin{equation}
\label{eq:def_ratio}
R(\alpha,t)=(1-\delta_R)\{(1-e^{-dt^2})R_e +e^{-dt^2}R_0\}+\delta_R R_0,\quad d>0.
\end{equation}
Here, $R_0$ and $R_e$ are defined as
%%---Def. of Components in Ratio Function---%%
\begin{equation}
\label{eq:def_ratio_comp}
R_0(\alpha,t)=\frac{1}{\pi}, \quad R_e(\alpha,t)=\frac{GL(s(\alpha),t)^{-1}}{\int_0^\pi GL(s(\alpha'),t)^{-1} d\alpha'},
\end{equation}
and these definitions trivially guarantee that 
 %%---Def. of s_\alpha with Ratio Function.---%%
\begin{equation}
\label{eq:ratio_total}
\int_0^{\pi} R_0d\alpha= \int_0^{\pi} R_ed\alpha= \int_0^{\pi} Rd\alpha=1.
\end{equation} 
Moreover, it is easy to show that the definition (\ref{eq:def_ratio}) satisfies the relations
%%---Value of Ratio Function at t=0---%%
\begin{align}
\label{eq:def_ratio_t0}
R(\alpha,0)=R_0,\quad \frac{\partial R}{\partial t}(\alpha,0)=0.
\end{align}
Hence, at $t=0$, the parametrization is uniform and the tangential velocity $V$ is given by the formula (\ref{eq:exp_uniV}). At $t\neq 0$, however, the equation (\ref{eq:inteq_V}) is no longer solvable explicitly for $V$. To circumvent this problem, NS04 suggests the backward difference approximation
%%---Approx. of dR/dt---%%
\begin{equation}
\label{eq:def_ratio_approx}
\frac{\partial R}{\partial t}(\alpha,t)\approx \frac{R(\alpha,t) - R(\alpha,t-\tau)}{\tau},
\end{equation}
for some $\tau>0$. This formula allows to find $V$ in the same way as for the time-independent $R$, although errors of $\mathcal{O}(\tau)$ are introduced to the relative spacing $s_\alpha$. We specify the values of $\tau$ in Section \ref{sec4}, and a validation of Eq.(\ref{eq:def_ratio_approx}) is given in Section \ref{sec5}.

%%%%%%%%%%%%%%%%%%%%%%%%%%%%%%%%%%%
%%%%%%%%%%%%%%%%%%%%%%%%%%%%%%%%%%%
%%%%%%      NUMERICAL METHODS         %%%%%%%%%%%%
%%%%%%%%%%%%%%%%%%%%%%%%%%%%%%%%%%%
%%%%%%%%%%%%%%%%%%%%%%%%%%%%%%%%%%%

\section{Numerical methods}\label{sec4}
We numerically solve the initial value problem of Eqs.(\ref{eq:dt_gamma})(\ref{eq:ev_theta})(\ref{eq:evol_sa}) with the normal velocity $U=\mathbf{W}\cdot \mathbf{n}$ and the tangential velocity $V$ specified in Section \ref{sec3}. Our numerical method consists mainly of spatial discretization, evaluations of the guideline function $GL$, and temporal discretization compatible with the approximation (\ref{eq:def_ratio_approx}).
\subsection{Spatial discretization}
In our computations, the dynamical variables $(\theta, s_\alpha, \gamma)$ are sampled at $N$ points equally spaced in the extended parameter space $[0,2\pi]$:
%%---Discretization of [0,2\pi] ---%%
\begin{equation}
\label{eq:def_grid}
\alpha_j = jh,\quad h=\frac{2\pi}{N}, \quad j=0,1,\ldots, N-1,
\end{equation}
and all the derivatives and anti-derivatives are evaluated using the standard fast Fourier transform (FFT). For the principal-values (\ref{eq:def_axbiot_r}) and (\ref{eq:def_axbiot_z}), we employ the third-order modified dBM approximation \cite{Nitsche1999} with 8193 quadrature nodes. The values of $(r, z, \gamma)$ between the grid points (\ref{eq:def_grid}) are obtained by the Fourier interpolation after the reconstruction from $(\theta,s_\alpha)$. In our code, evaluations of the integrands in (\ref{eq:def_axbiot_r}) and (\ref{eq:def_axbiot_z}), including the iterative algorithms for the complete elliptic integrals $K$ and $E$ \cite{AbSt1965}, are parallelized on a single NVIDIA Tesla P100, a Graphical Processing Unit (GPU) for scientific computing in the double-precision arithmetics. Also,  the direct summation after the integrand evaluations is performed with the warp-shuffle reduction algorithm \cite{ChGrMc2014} on the same GPU. Here, we are aware that higher-order quadrature rules for the axisymmetric Biot-Savart integrals have been suggested \cite{NiBa1998,Nitsche2001} and successfully implemented by Nie01 and NS04 for the case with surface tension. However, those methods are originally designed for studies on an axisymmetric analogue of the Moore's singularity, which requires to approximate the principal values beyond the double-precision accuracy. In the present work, it is found that the method described above yields satisfactory results in validating our mesh refinement scheme and is therefore useful in practice.  \par
\subsection{Evaluating $GL$}
Next, we apply the trapezoidal rule to the Fourier coefficients (\ref{eq:def_ft_arc}) after the change of variables (\ref{eq:changevar}):
%%---DFT in Arclength---%%
\begin{equation}
\label{eq:def_dft_arc}
\hat{f}(k)\approx \frac{h}{L_p} \sum_{j=0}^{N-1} (f(s(\alpha_j))s_\alpha (\alpha_j)) e^{-2\pi i \frac{s(\alpha_j)}{L_p}}.
\end{equation}
Unfortunately, this discrete sum cannot be directly computed using the standard FFT, because the nodes $s(\alpha_j)$ are redistributed non-uniformly in $[0,L_p]$ by mesh refinement. Instead, we can use the non-uniform fast Fourier transform (NUFFT), which is an $\mathcal{O}(N\log N)$ algorithm approximating the following types of sums to a prescribed relative accuracy $\epsilon_\text{rel}$: 
%%---Types of Nonuniform FFT ---%%
\begin{alignat}{2}
\label{eq:def_nufft_type1}&\hat{f}_k=\sum_{j=0}^{N-1}f_je^{-ik x_j}, \quad &k=-\frac{M}{2},\ldots, \frac{M}{2}-1\quad &\mbox{(Type-1)}, \\
\label{eq:def_nufft_type2}&f_j=\sum_{k=-\frac{M}{2}}^{\frac{M}{2}-1}{\hat{f}_k}e^{ik x_j}, \quad &j=0,1,\ldots, N-1\quad &\mbox{(Type-2)}.
\end{alignat}
It is easy to see that the formula (\ref{eq:def_dft_arc}) naturally fits in the framework of the Type-1 NUFFT (\ref{eq:def_nufft_type1}). On the other hand, once the Fourier coefficients of $GL$ are obtained, it can be interpolated back to the non-uniform grid $s(\alpha_j)$ using the Type-2 NUFFT (\ref{eq:def_nufft_type2}). Here, we note that there are several implementations of the NUFFT algorithms for general purposes, including CMCL NUFFT \cite{GrLe2004}, NFFT3 \cite{KeKuPo2009}, and FINUFFT \cite{BaMaKl2019}. In our code, we employ FINUFFT with a few cores of Intel Xeon Gold 6136, for, to our knowledge, it is currently the fastest library parallelized with OpenMP. Also, it should be noted that the Type-1 NUFFT itself is not related to any quadrature rules and simply approximates sums of a particular form. In our scheme, the coefficients $s_\alpha h$ serve as quadrature weights compatible with the non-uniform nodes $s(\alpha_j)$ that evolve in time, which automatically solves this issue. \par
Our practical procedure for evaluating $GL$ is as follows. Firstly, we decide the largest wavenumber $k_\text{max}$ of the Fourier coefficients $\hat{\kappa}_z(k)$ computed for the mesh refinement purpose. This is equivalent to choosing a target spatial resolution in terms of the uniform grid that we aim to achieve with a non-uniform grid realized by our mesh refinement. Therefore, the number $k_\text{max}$ should be related to the parameter $\delta_R$ in (\ref{eq:def_ratio}) so that the target resolution will always be no finer than the best possible local resolution determined by $\delta_R$. Secondly, if we assume that our scheme is working ideally, the locally refined grid resolves $\kappa_z$ successfully, whereas complex exponentials, which oscillate uniformly in $[0,L_p]$, can be underresolved for sufficiently large $k$. To resolve both $\kappa_z$ and those complex exponentials, the functions $\kappa_z$, $s_\alpha$, and $s$ are interpolated by the truncated Fourier serires to an upsampled grid with $N_\text{up}$ points in $[0,2\pi]$, and then the coefficients $\hat{\kappa}_z(k)$ are computed using the Type-1 NUFFT. Thirdly, once the coefficients $\hat{\kappa}_z(k)$ are obtained for $|k| \leq k_\text{max}$, the curvature $\kappa_z$ as a function of the arclength and its Hilbert transform are computed by inverting the corresponding Fourier coefficients with the standard FFT, and the analytic envelope $E_r[\kappa_z]$ is evaluated in $[0,L_p]$. Lastly, we transform $E_r[\kappa_z]$ back to the Fourier space, apply the Gaussian filter, and finally interpolate the resulting $GL$ to the non-uniform grid $s(\alpha_j)$ using the Type-2 NUFFT. Since each step involves $\mathcal{O}(N)$ and $\mathcal{O}(N\log N)$ operations only, additional costs for our mesh refinement strategy are $\mathcal{O}(N\log N)$.
\subsection{Temporal discretization and filtering}
We evolve the dynamical variables $(\theta,s_\alpha, \gamma)$ using the classical fourth-order Runge-Kutta method. With surface tension, it is well-known that high-order terms in $\kappa_\alpha$ introduce stiffness into the governing equations and the stepsize $\Delta t$ is required to satisfy a severe stability constraint. In the case of axisymmetric vortex sheets, this constraint is known to be $\Delta t \leq C  \Delta s_\text{min} ^{3/2}$, where $\Delta s_\text{min}$ is the minimum distance between computational points in the symmetry plane. In this paper, we follow NS04 and try to satisfy the condition with $C=2.5$ as long as possible. For a general time-independent $R$, a semi-implicit scheme that remove the stiffness has been developed by HLS97, and an extension of this method to more general cases is a direction of our future works. \par
For the approximation (\ref{eq:def_ratio_approx}), we choose $\tau=c_i \Delta t$ at the $i$-th Runge-Kutta stage with $c_1= 1$, $c_2=c_3=1/2$, and $c_4=1$. Combining these coefficients with the property (\ref{eq:def_ratio_t0}), we can perform time-stepping by the classical Runge-Kutta method without any interpolation technique in the temporal direction. It is not clear that the coefficients used in NS04 are the same as our choice, and there may be some alternative values. A drawback of this procedure is that it amounts to introducing a highly oscillatory delay $\tau(t)$ inevitably associated with $\Delta t$ and may affect the convergence of our temporal discretization. This issue is carefully discussed in Section \ref{sec5}. \par
Also, the spectral filter introduced by Krasny \cite{Krasny1986} is applied to numerical solutions at every time step. This technique removes all Fourier coefficients whose magnitudes are below a cut-off level $\epsilon_K$, and it is originally used to control unphysical growths of round-off errors due to the Kelvin-Helmholtz instability in the zero surface tension case. In the present case, another source of noises is the principal curvature $\kappa_r$, which is a quotient of two functions that vanish on the axis of symmetry. Since Eq.(\ref{eq:dt_gamma}) involves the derivative of this function, large noises can be introduced to numerical solutions and grow under, for example, the Kelvin-Helmholtz instability. The Krasny's filter is also expected to control the latter problem, as demonstrated in NS04.

%%%%%%%%%%%%%%%%%%%%%%%%%%%%%%%%%%%
%%%%%%%%%%%%%%%%%%%%%%%%%%%%%%%%%%%
%%%%%%      NUMERICAL RESULTS          %%%%%%%%%%%%
%%%%%%%%%%%%%%%%%%%%%%%%%%%%%%%%%%%
%%%%%%%%%%%%%%%%%%%%%%%%%%%%%%%%%%%

\section{Numerical results}\label{sec5}
Now we show some numerical results to validate the numerical method described in Section \ref{sec4}. In the following computations, the parameters for our mesh refinement scheme are given by
\begin{equation}
\label{eq:def_ref_param}
a=20,\quad d=5,\quad \delta_R=\frac{1}{8},\quad  k_\text{max} = 8 \biggl(\frac{N}{2}\biggr), \quad N_\text{up}=32N,
\end{equation}
with the relative accuracy $\epsilon_\text{rel}=10^{-15}$ for the NUFFT algorithms. The other parameters $N$, $\Delta t$ and $\epsilon_K$ are specified in each subsection.
%
%
%%%%%%      INITIAL PARAMETRIZATION         %%%%%%%%%%%%%%%%%
%
%
\subsection{Construction of uniform parametrization}
    As a limitation of the definition (\ref{eq:def_ratio}), it is required to construct an initial condition in the uniform parametrization before any time-stepping is performed. There are at least two ways to complete this step. The first is to determine the total length $L$ (or equivalently the relative spacing $s_{\alpha}$ constant in $\alpha$) and the functions $\theta$ and $\gamma$ that are assumed to be in the uniform representation. This approach is cleverly used by Almgren \cite{Almgren1996} to design initial shapes of Hele-Shaw bubbles that can lead to finite-time singularities. The second is to find the uniform parametrization of a given initial condition using some numerical techniques. Fortunately, the latter approach can be implemented by an immediate application of the ideas described in Section \ref{sec4}. Remember that the formula (\ref{eq:def_dft_arc}) is the trapezoidal rule applied to a smooth periodic function of the variable $\alpha$, which is known to be spectrally accurate. On the other hand, the arclength parametrization is a special case of the uniform ones whose domain is $[0,L_p]$. Hence, by inverting the approximate coefficients (\ref{eq:def_dft_arc}) via the standard FFT, we expect to obtain a function represented in the uniform parametrization with spectral accuracy. To our knowledge, this procedure is distinctly different from the existing methods designed for the same purpose. For instance, Baker and Shelley \cite{BaSh1990} suggest an iterative scheme based on the Newton's method, which searches for a set of values $\alpha_j$ corresponding to $s(\alpha_j)$ equally spaced in $[0,L_p]$. In another direction, Mikula and \v{S}ev\v{c}ovi\v{c} \cite{MiSe2004} propose the concept called {\it asymptotically uniform parametrization}, which can define a ``time-evolving" method that fixes $U=0$ and evolves a plane curve solely by the tangential velocity $V$ until the relative spacing $s_\alpha$ converges to a constant. Our method needs neither the Newton's method nor temporal discretization, and it is completely non-iterative. \par 
To justify the argument above, our initialization algorithm is tested on the following example in polar coordinates $(\eta,\phi)$:
%
%%---Periodic Solution in Linear Problem---%%
\begin{equation}
\label{eq:linear_period}
\eta(\phi)=1+\epsilon_P P_2(\cos\phi), \;\; \phi \in [0,2\pi].
\end{equation} 
Here, $\phi$ is the polar angle, $\eta$ is the radius as a function of $\phi$, and $P_2$ is the Legendre polynomial of order two. This special choice, which is presented in \cite{Lamb1932}, corresponds to the zero-velocity state of a periodic solution to the linearized problem in \cite{Rayleigh1879}. The representation in the polar coordinates allows to compute the ``inverse" mapping $(r,z)\mapsto \phi$ and check the accuracy of numerical results against the analytical data (\ref{eq:linear_period}), although our method is also applicable to smooth closed curves given in Cartesian coordinates.\par
\begin{figure}[t]
\begin{center}
\includegraphics[width=\linewidth,trim=0 0 0 0]{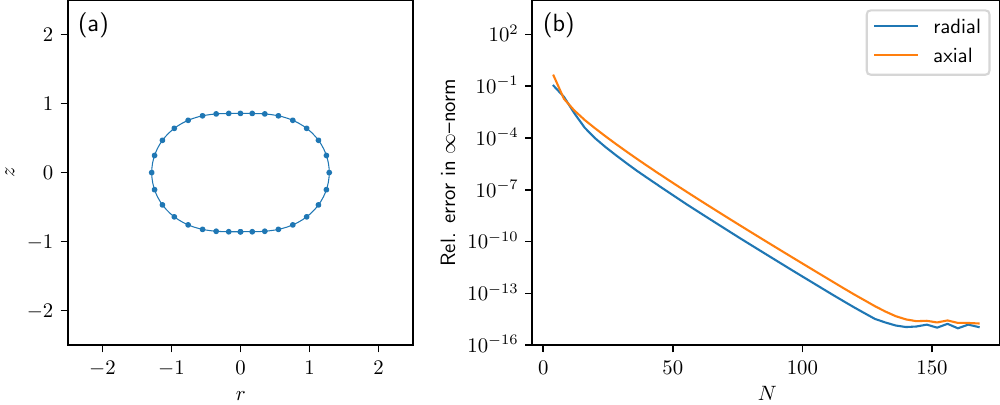}
\end{center}
\caption{\label{fig:conv_init} (a) Interface position of test problem (\ref{eq:linear_period}) in $r$--$z$ plane, (b) exponential convergence of initialization algorithm based on inverting approximate coefficients (\ref{eq:def_dft_arc}).}
\end{figure}

Following HLS94, we first compute the curvature $\kappa_z$ and the relative spacing $s_\alpha$ from (\ref{eq:linear_period}), and find the representation of $\kappa_z$ in the uniform parametrization by our initialization algorithm. Then, with $s_\alpha=L_p/2\pi$, the derivative $\theta_\alpha=s_\alpha\kappa_z$ is integrated to obtain $\theta$, which is followed by the reconstruction of $(r,z)$. In Fig.\ref{fig:conv_init}, the plane curve (\ref{eq:linear_period}) and the convergence of the scheme are shown for $\epsilon_P=2/7$. 
%
%%---Fig. Spectral Accuracy of the Construction---%%
As shown in Fig.\ref{fig:conv_init}(a), it has the variable radius $\eta$ and therefore its parametrization is non-uniform. Next, in Fig.\ref{fig:conv_init}(b), numerical errors in $r$ and $z$ are plotted versus $N$ in terms of the $\infty$--norm divided by the maximal value of each function. The exponential convergence of the trapezoidal rule is clear, and each graph reaches a level slightly above $10^{-15}$ for $N\geq 150$. In fact, we do not perform a construction of the uniform parametrization for the initialization purpose, because the conditions (\ref{eq:init_nie2}) and (\ref{eq:init_nie1}) are uniform by definition. Instead, our algorithm is effectively used in the convergence study of the time-stepping with the approximation (\ref{eq:def_ratio_approx}), which requires to remove effects of mesh refinement for comparisons between two solutions with different $\Delta t$.
%
%
%
%%%%%%      CAPILLARY PINCHOFF         %%%%%%%%%%%%%%%%%
%
%
\subsection{Capillary pinch-off}
 Our numerical method is first tested on the pinch-off problem (\ref{eq:init_nie2}) with $N=2048$, $\Delta t=5\times 10^{-5}$, and $\epsilon_K = 10^{-11}$.  
 %
%%---Fig. Snapshots of Capillary Pinchoff---%%
\begin{figure}[t]
\begin{center}
\includegraphics[width=\linewidth,trim=0 0 0 0]{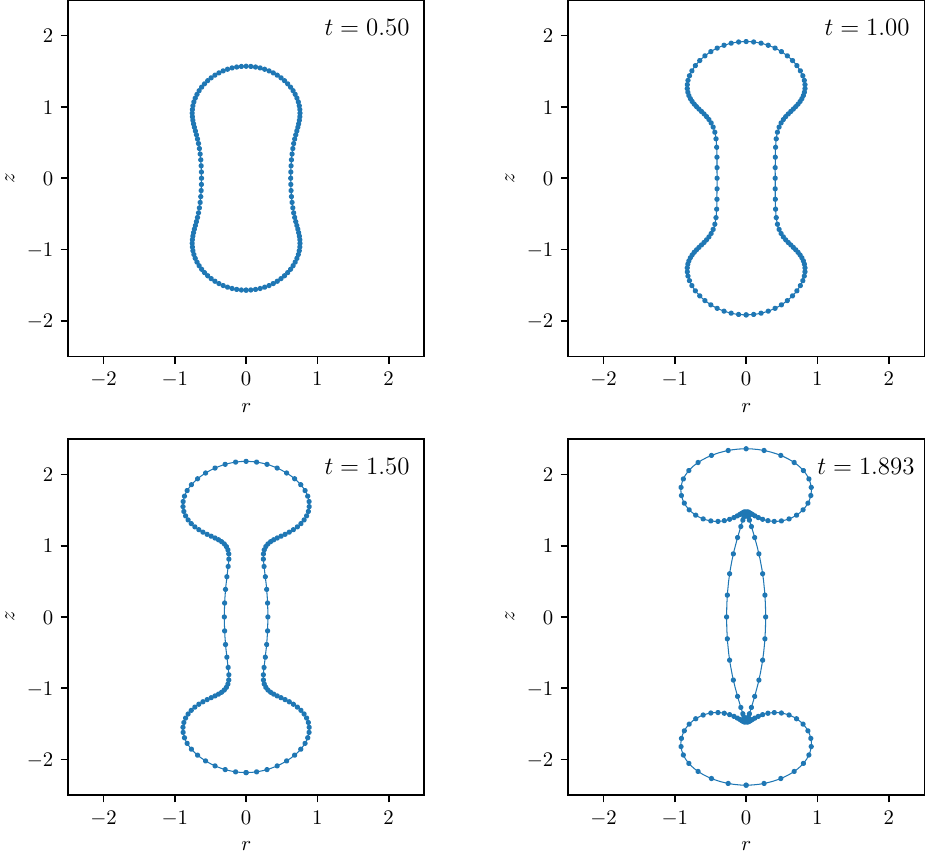}
\end{center}
\caption{\label{fig:snaps_pinchoff} Snapshots of interface position at indicated times for pinch-off problem (\ref{eq:init_nie2}).}
\end{figure}
 Fig.\ref{fig:snaps_pinchoff} plots the numerical solutions in the $r$--$z$ plane at indicated times, including the symmetric image in $r<0$. To clearly show arrangements of computational points, in each figure the solid line is drawn at the the full resolution and every 16th point is plotted again as a dot.  At $t=0.50$, the droplet starting from a sphere is axially deformed and acquires a capsule-like shape. Since it is still relatively close to a sphere, the effect of mesh refinement is only slightly visible as sparseness near $z=0$. At $t=1.00$, the capsule-like droplet is further elongated along the $z$--axis, and now it has developed two round ends connected by a thick rod. At this moment, the rod in the middle is almost straight near $z=0$, while geometry is relatively simple around the end-points $\mathbf{X}(0,t)$ and $\mathbf{X}(\pi,t)$. For these simplicities, computational points are beginning to accumulate near the junctions between the rod and the round ends, as clearly seen in the corresponding figure. Almost the same argument applies to the solution at $t=1.50$. Here, the radius of the rod is no longer uniform, and two isolated points of the minimum neck radius are formed near the junctions, implying where pinch-off is about to occur. At $t=1.893$, which is just before noises become visible in the numerical solution, the minimum neck radius has rapidly decreased and the droplet is forming two end-droplets at the top and bottom along with a sharply elongated satellite droplet in between. Also, as discussed in Section \ref{sec2}, overturning is found in each pinch-off region, and it develops the so-called cone-crater structure carefully studied by NS04. For reference, this GPU-accelerated simulation took approximately 1.25 hours with a single CPU core and one hour with two CPU cores. We could not find significant speed gain from further increasing the numbers of cores, which is presumably because the orders of $k_\text{max}$ and $N_\text{up}$ were relatively small compared to the main target of the parallelized NUFFT algorithms (see \cite{BaMaKl2019} for details).  The same simulation was performed with the uniform parametrization and took roughly 10 minutes, although the numerical solution at around $t=1.89$ was contaminated by noises at the highest wavenumber.\par 
%
%%---Fig. Snapshots of Capillary Pinch-off---%%
\begin{figure}[t]
\begin{center}
\includegraphics[width=\linewidth,trim=0 0 0 0]{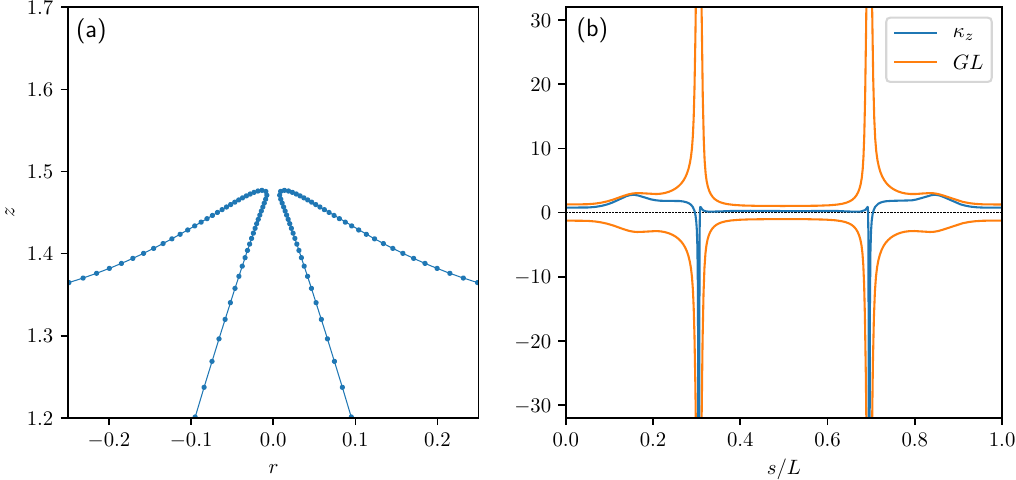}
\end{center}
\caption{\label{fig:closeup_pinchoff} Interface profile for pinch-off problem (\ref{eq:init_nie2}) at $t=1.893$: (a) closeup of pinch-off region with every 4th point plotted as dot, (b) plots of $\kappa_z$ and $GL$.}
\end{figure}
In Fig.\ref{fig:closeup_pinchoff}(a), we show a closeup of the upper pinch-off region at $t=1.893$. As easily seen in the figure, a small region containing the shrinking neck is densely resolved by our mesh refinement scheme, and the minimum value of  $s_\alpha$ is six times smaller than that of the uniform parametrization. This refinement is roughly 1.3 times coarser than that of HLS97, and two times than that of NS04. In Fig.\ref{fig:closeup_pinchoff}(b), the curvature $\kappa_z$ and the corresponding $GL$ are shown for the same $t$. As expected, the procedure outlined in Section \ref{sec4} successfully generates a simple and strictly positive $GL$ from the curvature $\kappa_z$ with two significantly sharp peaks. Here, note that the locality of mesh refinement can be improved by increasing the parameter $a$ in (\ref{eq:def_heat}) at the cost of the smoothness of $GL$. An ideal scheme should be capable of updating the value of $a$ adaptively in the course of simulations as high-curvature regions develop. Unfortunately, we have not found any indicator that can evolve $a(t)$ in a differentiable manner without any a priori information, and the fixed parameters in (\ref{eq:def_ref_param}) have been obtained by trial and error. However, once we know a set of parameters that works properly, it can often be reused for other initial conditions, as we will see later in the second example (\ref{eq:init_nie1}). \par
%
%%---Fig. Fitting Self-Similarity ---%%
\begin{figure}[t]
\begin{center}
\includegraphics[width=\linewidth,trim=0 0 0 0]{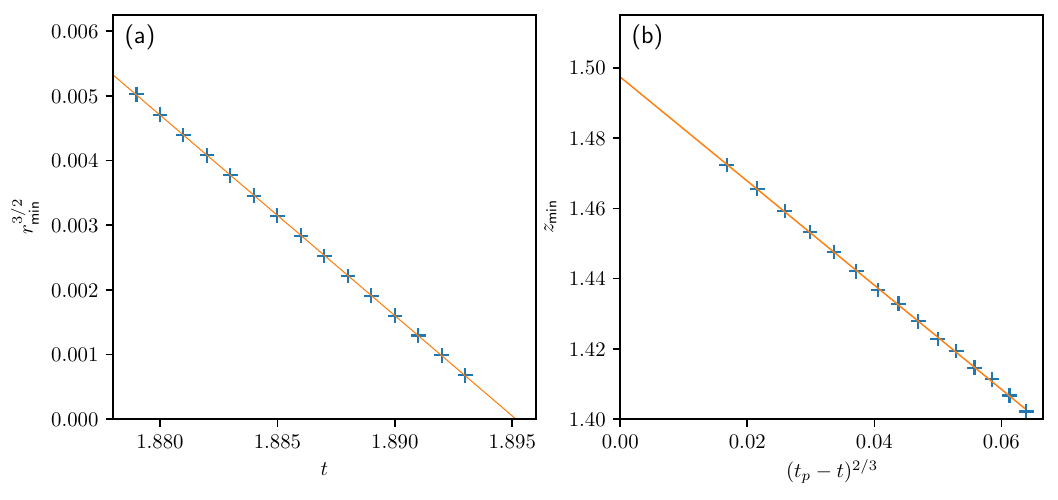}
\end{center}
\caption{\label{fig:scaling_rz}Confirmation of scaling law (\ref{eq:scaling_nie2}) for pinch-off problem (\ref{eq:init_nie2}): (a) $r_\text{min}^{3/2}$ vs $t$, (b) $z_\text{min}$ vs $(t_p-t)^{2/3}$.}
\end{figure}
Next, we confirm the scaling law (\ref{eq:scaling_nie2}) using our numerical results. The procedure is essentially the same as in NS04. Firstly, in Fig.\ref{fig:scaling_rz}(a), we plot $r_\text{min}^{{3/2}}$ versus $t$ and fit a straight line to the data points. As one can easily see, the fitted line clearly explains the linear behavior of $r_\text{min}^{{3}/{2}}$, which suggests that the solution is in the self-similar regime. Secondly, by extrapolating the line to $r^{{3}/{2}}_\text{min}=0$, we obtain an estimate $t_p\approx1.8951$, which is slightly smaller than $t_p\approx1.89523$ obtained in NS04. Thirdly, the coordinate $z_\text{min}$ is plotted in Fig.\ref{fig:scaling_rz}(b) as a function of $(t_p-t)^{{2}/{3}}$ using the value of $t_p$ estimated at the previous step, and then another straight line is fitted to the data points. Although it is not as clear as for $r^{{3}/{2}}_\text{min}$, it is still reasonable to conclude that the results are well explained by the fitted line. Again, by extending the line to $(t_p-t)^{{2}/{3}}=0$, we obtain an estimate $z_p\approx 1.4973$, while that by NS04 is $z_p\approx 1.49839$. From these results, we believe that  our simulation with $N=2048$ resolves the singularity formation slightly better than one by NS04 with the uniform parametrization and $N\approx 8000$. However, it should be also pointed out that the deviations of our estimates from those in NS04 are arguably due to the fact that our mesh refinement scheme fails to continue simulations at further smaller length-scales, where simple universal scalings such as (\ref{eq:scaling_nie2}) can be distinguished from more complicated situations including bubble collapse dominated by inertia \cite{EgFoLeSn2007} and droplet formations affected by weak viscosity \cite{DeHeHaVeRoKeEgBo2018}. The failure of our simulation is attributed mainly to the lack of an adaptive stepsize control, which is required to handle large interfacial velocity near singularity formations as well as the stiffness $\Delta t \leq C  \Delta s_\text{min} ^{3/2}$ that becomes more severe as the spatial resolution is improved. This drawback is due to the introduction of the time delay $\tau$ in (\ref{eq:def_ratio_approx}) linked to the stepsize and the classical Runge-Kutta method, which prevents us from updating $\Delta t$ in a meaningful manner. Regarding this point, a promising strategy may be to use some interpolation technique in the temporal direction so that the delay $\tau$ can be untied from $\Delta t$ and fixed at a sufficiently small value. Such flexibility may also allow to use more stable time-integration schemes than the classical Runge-Kutta method. We hope to address this issue in future investigations.  \par
%
%%---Fig. Geometric Convergence of Runge-Kutta method---%%
\begin{figure}[t]
\begin{center}
\includegraphics[width=\linewidth,trim=0 0 0 0]{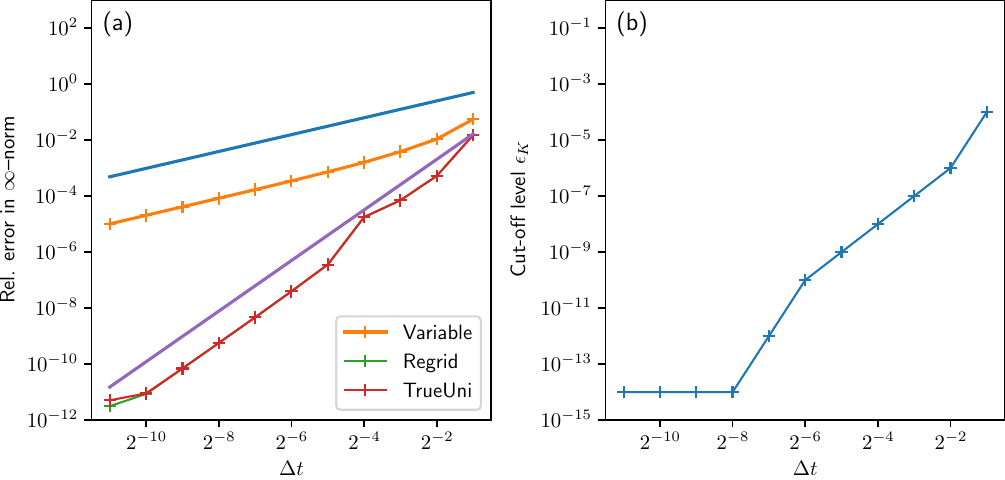}
\end{center}
\caption{\label{fig:conv_geom} Convergence of time-stepping with approximation (\ref{eq:def_ratio_approx}): (a) plots of errors measured as relative $\infty$--norm, (b) correspondence between stepsize $\Delta t$ and filter level $\epsilon_K$.}
\end{figure}
%
%
%
%
%%%%%%      CONVERGENCE STUDY         %%%%%%%%%%%%%%%%%
%
%
\subsection{Convergence study}
NS04 justifies the controversial approximation (\ref{eq:def_ratio_approx}) by showing that changing the resolution $N$ or the parameter $\delta_\text{max}$ in (\ref{eq:def_NiSt_guide}) does not affect the results in their figures similar to Fig.(\ref{fig:scaling_rz}). In this paper, we aim to verify in a more direct way that the time-stepping involving (\ref{eq:def_ratio_approx}) actually converges faster than $\mathcal{O}(\Delta t)$.  Before doing any special trick, we first perform a standard convergence study with the pinch-off problem (\ref{eq:init_nie2}) and the time $t=0.50$ (see Fig.\ref{fig:snaps_pinchoff}), in which solutions with various $\Delta t$ are compared in $[0,\pi]$ against an ``exact" solution obtained with a sufficiently small stepsize. Here, to weaken the stiffness $\Delta t \leq C  \Delta s_\text{min} ^{3/2}$, the number $N$ is reduced to 256, and the ``exact" solution is computed with $\Delta t=10^{-5}$ and $\epsilon_K=10^{-14}$. For large $\Delta t$ that still violate the constraint, the instability typically appears as a rapid growth in the Fourier coefficients of numerical solutions at high wavenumbers. This fact is numerically demonstrated by HLS94 for two-dimensional vortex sheets with surface tension, and the same phenomenon should be seen in the axisymmetric case. To continue simulations for such large $\Delta t$, we control the growth of the Fourier coefficients by adjusting the filter level $\epsilon_K$. In Fig.\ref{fig:conv_geom}(a), we plot errors in the mapping $\mathbf{X}$ measured as the $\infty$--norm divided by the maximum value of $|\mathbf{X}(\alpha_j,t)|$, and the correspondence between $\Delta t$ and $\epsilon_K$ is shown in Fig.\ref{fig:conv_geom}(b). As seen in Fig.\ref{fig:conv_geom}(a), the error curve for $\mathbf{X}$ (denoted by ``Variable" in the legend) is parallel to a linear function (blue line) and the convergence rate is apparently $\mathcal{O}(\Delta t)$. From this result, one may be tempted to immediately conclude that the order of the classical Runge-Kutta method is completely lost. However, noting that a change in the stepsize $\Delta t$ can alter the relative spacing $s_\alpha$ linearly through the approximation (\ref{eq:def_ratio_approx}) with $\tau=c_i\Delta t$, we still hesitate to eliminate the possibility that the true convergence rate may be hidden behind the $\mathcal{O}(\tau)$ behavior of (\ref{eq:def_ratio_approx}).\par
From this point of view, we repeat the same convergence study after removing the effects of our mesh refinement scheme. That is, the computed solution for each $\Delta t$ is transformed to its arclength parametrization by the initialization technique based on inverting (\ref{eq:def_dft_arc}), and then errors are measured again as the relative $\infty$--norm in the same manner as above. This convergence study is  geometric in the sense that the difference between two plane curves is considered in a canonical form, which is the arclength parametrization in the present case. Again in Fig.\ref{fig:conv_geom}(a), we show some error curves obtained after removing potential artifacts due to the approximation (\ref{eq:def_ratio_approx}). Firstly, the error curve computed from the regridded solutions (denoted by ``Regrid") is plotted with a theoretical $\mathcal{O}(\Delta t^3)$ curve (purple curve).  As opposed to the result above, the convergence rate here is clearly superlinear and estimated to be at least $\mathcal{O}(\Delta t^3)$, although some order reduction may be occurring presumably due to the fast oscillation of $\tau=c_i\Delta t$. Also, to check that the numerical solutions with our mesh refinement converge to the same curve as with the uniform parametrization, we plot the error curve (denoted by ``TrueUni") obtained from comparisons between the regridded solutions and an ``exact" solution with the tangential velocity (\ref{eq:exp_uniV}). The two dotted curves are indistinguishable except for the left-most data point, which strongly suggests that the numerical method described in Section \ref{sec4} converges to a correct solution faster than the order of the approximation (\ref{eq:def_ratio_approx}). 
%
%%%%%%      INERTIA BREAKUP         %%%%%%%%%%%%%%%%%
%
%
\subsection{Bag breakup}
\begin{figure}[t]
\begin{center}
\includegraphics[width=\linewidth,trim=0 0 0 0]{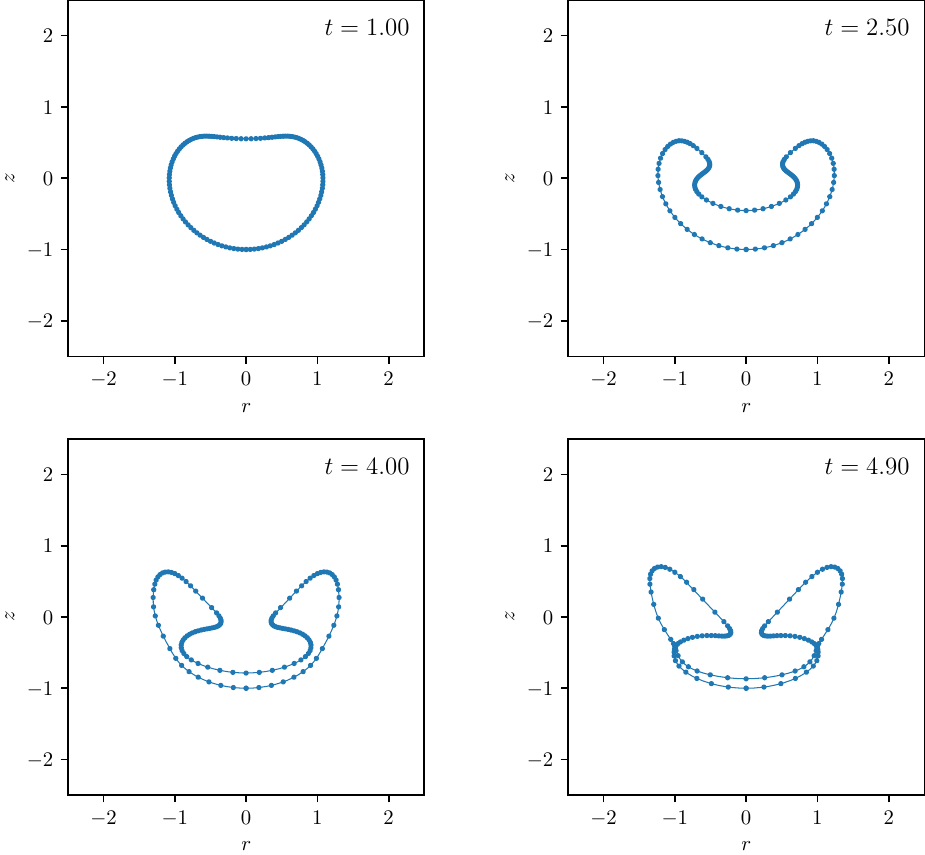}
\end{center}
\caption{\label{fig:snaps_inertia} Snapshots of interface position at indicated times for axisymmetric bag breakup (\ref{eq:init_nie1}).}
\end{figure}

Next, the performance of our numerical method is further examined in its application to the axisymmetric bag breakup (\ref{eq:init_nie1}) with the same $N$, $\Delta t$, and $\epsilon_K$ as for the pinch-off problem. Again, Fig.\ref{fig:snaps_inertia} plots the numerical solutions in the $r$--$z$ plane (including $r<0$) at indicated times, and every 16th point is plotted as a dot after drawing the solid line at the the full resolution. 
%
%%---Fig. Snapshots of Inertia Breakup---%%
Also, since we do not evolve the end-point $\mathbf{X}(0,t)$ in this simulation, the interface position in each figure is reconstructed with $\mathbf{X}(0,t)=(0,-1)$. At $t=1.00$, the droplet, which is initially a sphere, has developed a small indentation at the top and the inner region is no longer convex. The effect of mesh refinement is not evident here because no complex geometry is found in the interface.  At $t=2.50$, as the indentation further develops, the droplet begins to acquire a bowl-like shape whose edge is bent toward the axis of symmetry. Now, computational points weakly cluster near the two round corners of the edge, and consequently some sparseness is found away from those small regions.  At $t=4.00$, a large volume of the inner fluid has flowed into the edge region and the bottom of the bowl has become thinner.  At $t=4.90$, which is just before the numerical solution is contaminated by noises, the fluid interface is about to collide itself at the point connecting the edge and the bottom. As one can see, most of computational points are now devoted to resolving the two corners of the edge and the small region around the self-intersection.\par
\begin{figure}[t]
\begin{center}
\includegraphics[width=\linewidth,trim=0 0 0 0]{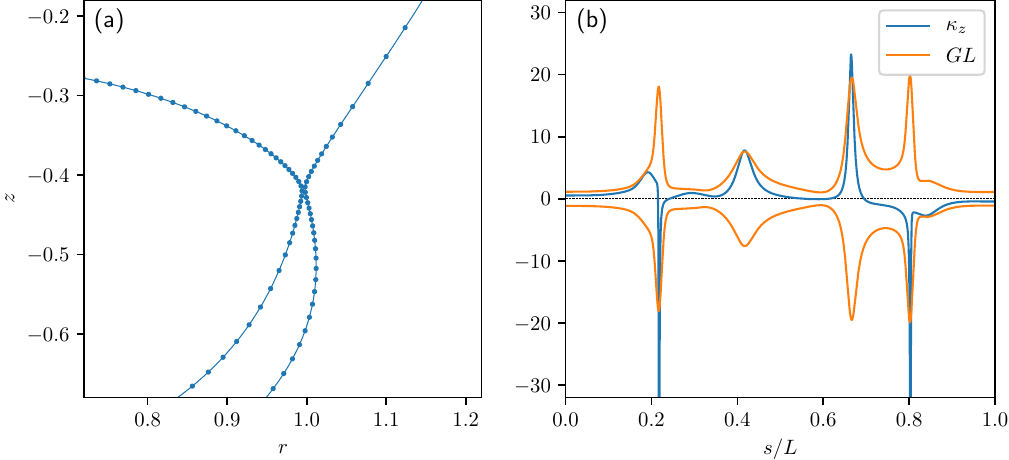}
\end{center}
\caption{\label{fig:closeup_inertia} Interface profile for axisymmetric bag breakup (\ref{eq:init_nie1}) at $t=4.90$: (a) closeup of self-intersection region with every 4th point plotted as dot, (b) plots of $\kappa_z$ and $GL$.}
\end{figure}

 
The topological change shown in Fig.\ref{fig:snaps_inertia} should be considered qualitatively different from that in Fig.\ref{fig:snaps_pinchoff}. In the case of bag breakup, the droplet is about to break into one torus-like component and one simply-connected component, whereas the symmetric pinch-off seems to form three simply-connected ones. In Fig.\ref{fig:closeup_inertia}(a), we show a closeup of the region near the point of the self-intersection at $t=4.90$.
%
%%---Fig. Snapshots of Inertia Breakup---%%
As opposed to capillary pinch-off, overturning or a cone-crater structure is not observed in this case, and the fluid interface has two corners facing each other. This structure is more similar to those found in, for example, HLS97 and Burton and Taborek \cite{BuTa2007B}. For the same $t$, the curvature $\kappa_z$ and the corresponding $GL$ are shown in Fig.\ref{fig:closeup_inertia}(b). In this case, the solution forms two high-curvature regions away from the points of the self-intersection, which precludes applications of mesh refinement based on the prescribed guideline (\ref{eq:def_NiSt_guide}). As expected, our method succeeds in generating a simple and strictly positive $GL$ from the curvature $\kappa_z$ with four peaks with different magnitudes. Here, the relative spacing $s_\alpha$ is roughly two times smaller for the smallest peak than that of the uniform parametrization, and 3.5 times smaller for the others three. On the other hand, a drawback of mesh refinement based on the analytic envelope is found in the interval between the two right-most peaks in Fig.\ref{fig:closeup_inertia}(b). In the function $GL$, the two adjacent peaks are ``bridged", which is again presumably due to the low-frequency nature of the operator $E_r$,  and our mesh refinement scheme overestimates the complexity of a region where the curvature $\kappa_z$ is small and smooth (see Fig.\ref{fig:snaps_inertia}). This phenomenon is not due to the effect of the Gaussian filter. In fact, this bridged function is found ``before" applying any smoothing filter, and therefore inherent to the use of the analytic envelope.   

%%%%%%%%%%%%%%%%%%%%%%%%%%%%%%%%%%%
%%%%%%%%%%%%%%%%%%%%%%%%%%%%%%%%%%%
%%%%%%      CONCLUSION                        %%%%%%%%%%%%
%%%%%%%%%%%%%%%%%%%%%%%%%%%%%%%%%%%
%%%%%%%%%%%%%%%%%%%%%%%%%%%%%%%%%%%
\section{Conclusion}\label{sec6}
In this paper, we have proposed a new mesh refinement scheme for the vortex sheet formulation of inviscid droplet dynamics with axial symmetry. The method is based on signal processing using the Fourier coefficients of the curvature in the $r$--$z$ plane as a function of the arclength. The combination of the analytic envelope and the subsequent Gaussian filter enables to detect two kinds of singularity formations with far less human intervention than the existing methods, and additional computational costs for our scheme are $\mathcal{O}(N\log N)$ with the help of the non-uniform fast Fourier transform. We have also proposed a non-iterative construction of the uniform parametrization for closed plane curves. It is used to confirm, at least in a geometric sense, a convergence of the time-stepping procedure with the approximation (\ref{eq:def_ratio_approx}) introduced by NS04. \par
There are two important directions for future developments related to the issues discussed in Section \ref{sec5},. Firstly, a hanging bridge between two peaks of the guideline function $GL$, which is found in the case of bag breakup, results in some inefficiency of our mesh refinement scheme. This drawback motivates us to seek an alternative technique from signal processing that replaces the analytic envelope employed in this paper. Also, for the smoothing step, it may be interesting to switch from the heat kernel to other window functions such as the Kaiser‐Bessel family \cite{Kaiser1966}. Secondly, an adaptive stepsize control is desirable in simulations of vortex sheets with surface tension. In addition to small time-scales near self-intersections and the stiffness of the governing equations, HLS94 reports that the ``shadow" of the Moore's singularity and its saturation will lead to fast changes in solutions and can cause sudden drops in overall accuracy. Moreover, Ambrose \cite{Ambrose2009} points out using their model problem that a nearly-singular behavior unrelated to any self-intersection can appear even in the non-zero surface tension case. As discussed, one possibility is to interpolate numerical solutions in the temporal direction so that we can freely choose the delay $\tau$ regardless of the stepsize or the time-integration scheme. These projects are currently in progress.

%%%%%%%%%%%%%%%%%%%%%%%%%%%%%%%%%%%
%%%%%%%%%%%%%%%%%%%%%%%%%%%%%%%%%%%
%%%%%%      ACKNOWLEDGEMENT                        %%%%%%%%
%%%%%%%%%%%%%%%%%%%%%%%%%%%%%%%%%%%
%%%%%%%%%%%%%%%%%%%%%%%%%%%%%%%%%%%

\section*{Acknowledgement}
We would like to thank Jon Wilkening for providing us with helpful comments and computational resources at Berkeley. We are also grateful to Toshio Aoyagi for providing us with computational environments at Kyoto. This research was partly conducted in the semester program ``Singularities and Waves in Incompressible Fluids" at The Institute for Computational and Experimental Research in Mathematics (ICERM). We acknowledge the support from JSPS Overseas Challenge Program for Young Researchers (No. 201880131).

\bibliographystyle{unsrt}
\bibliography{refs_sp}

\end{document}